[1994/12/01]
\documentclass[reqno,12pt]{amsart}
\usepackage{bbm}
\def\1{\mathbbm{1}}

\usepackage{amsmath,amsthm}
\usepackage[mathscr]{euscript}
\usepackage{mathrsfs}
\usepackage{amsfonts}
\usepackage{amssymb}
\usepackage{latexsym}
\usepackage[pdftex,bookmarks,colorlinks,linkcolor=blue]{hyperref}
\usepackage{hyperref}
\newtheorem{defi}{Definition}[section]
\newtheorem{lem}{Lemma}[section]
\newtheorem{prop}{Proposition}[section]
\newtheorem{theorem}{Theorem}[section]
\theoremstyle{remark}
\newtheorem{remark}{Remark}[section]
\newtheorem{example}{Example}[section]

\numberwithin{equation}{section}

\newcommand{\T}{\mathbb{T}}
\newcommand{\R}{\mathbb{R}}
\newcommand{\C}{\mathbb{C}}
\newcommand{\N}{\mathbb{N}}
\newcommand{\F}{\mathbb{F}}

\newcommand{\D}{\mathbb{D}}

\def\al{\alpha}

\def\si{\sigma}

\def\la{\lambda}

\def\t{\tau}

\newcommand{\A}{\mathcal{A}}
\newcommand{\B}{\mathcal{B}}
\renewcommand{\L}{\mathcal{L}}

\newcommand{\bpm}{\begin{pmatrix}}
\newcommand{\epm}{\end{pmatrix}}

\newcommand{\brm}{\begin{remarque}}
\newcommand{\erm}{\end{remarque}}

\title{Analysis and control of integro-differential Volterra equations with delays }
\author{Y. El Kadiri, S. HADD and H. Bounit}
\address{Department of Mathematics, Faculty of Sciences, Ibn Zohr University, Hay Dakhla, BP8106, 80000--Agadir, Morocco; y.elkadiri@gmail.com, s.hadd@uiz.ac.ma, h.bounit@uiz.ac.ma}

\hypersetup{
  colorlinks=true,
  linkcolor=blue, 
  citecolor=red, 
 urlcolor=blue  } 
\begin{document}
\maketitle

\renewcommand{\sectionmark}[1]{}
\begin{abstract}
The purpose of this paper is to introduce a semigroup approach to linear  integro-differential systems with delays in state, control and observation parts.  On the one hand,  we use product spaces to reformulate state-delay integro-differential equations to a standard Cauchy problem and then use a perturbation technique (feedback) to prove the well-posendess of the problem, a new variation of constants formula for the solution as well as some spectral properties. On the other hand, we use the obtained results to prove that integro-differential systems with delays in state, control and observation parts form a subclass of distributed infinite-dimensional regular linear systems in the Salamon-Weiss sense.
 \end{abstract}

\section{Introduction}\label{sec:1}
Integro-gdifferential equations have attracted the attention of many researchers for many years and have been the subject of much interesting work, see e.g. \cite{Mil-73},\cite{Mil-75},\cite{Pruss93}, and the references therein. In such a class of evolution equations, the differential and integral operators can appear together at the same time.

In the first part of this work, we are concerned with introducing a unified semigroup approach to the following class of delay-integro-differential equations
\begin{equation}\label{problem}
\begin{cases}
\dot{x}(t)=Ax(t)+\displaystyle\int_0^ta(t-s)Ax(s) ds+Lx_t+f(t),\qquad t\geq 0\\
x(0)=x,\quad x_0=\varphi
\end{cases}
\end{equation}
where $A$ is a generator of a strongly continuous semigroup $(T(t))_{t\ge 0}$ on a Banach space $X$,
$x\in X,\; \varphi\in L^p([-r,0],X),\;a\in W^{1,p}(\R_+,\C)$ with $p\in (1, \infty)$ and $L$
is a Riemann-Stieljes integral of the form
\begin{align}\label{L-RS}
L\varphi:=\int^0_{-r}d\mu(\theta)\varphi(\theta),\qquad \varphi\in W^{1,p}([-r,0],X),
\end{align}
with  $\mu:[-r,0]\to\mathcal{L}(X)$ is a
function of bounded variation continuous in $0$ and with total variation
 $|\mu|$ (which is a positive Borel measure) satisfying $|\mu|([-\varepsilon,0])\to 0$ as $\varepsilon\to 0$.
 Here for each $t\ge 0$, the history function $x_t:[-r,0]\to X$ is defined by $x_t(\theta)=x(t+\theta)$ for $\theta\in [-r,0]$.

Note that a complete theory already exists for the problem \eqref{problem} in the case $a\equiv 0$ (i.e. the standard delay equations, see e.g., \cite{BatPia} and \cite{Hadd-SF}),  and in the case $L\equiv 0$ (called the free-delay integro-differential equation, see e.g. \cite{DGS-88}, \cite{DS-85}, \cite[Chap.7]{EnNag} and \cite{NS-93}). A common technique for solving the two cases mentioned above is using product spaces and matrix operators, see e.g. \cite{BatPia},\cite[Chap.7]{EnNag} for an elegant presentation.  We also mention that the control theory of free-delay integro-differential equations with unbounded control operator was first studied in \cite{Ju-00} and investigated by the authors in \cite{JP-04}, \cite{japar}, \cite{FadiliBounitb} and \cite{FadiliBounita} .



In the present work, we will also use matrix operators and feedback theory of infinite-dimensional linear systems \cite{WeiRegu} to study the well-posedness and develop a spectral theory for the  delay-integro-differential equation \eqref{problem}. In fact, we introduce the  product space
\begin{align*}
\mathcal{X}=X\times L^p(\R^+,X)\times L^p([-r,0],X).
\end{align*}
Using this space, we embed the problem \eqref{problem} in a large free-delay Cauchy problem in $\mathcal{X}$ of the form
\begin{align*}
\dot{z}(t)=\mathfrak{A}_L z(t),\quad z(0)=\left(\begin{smallmatrix}x\\ f\\\varphi\end{smallmatrix}\right),\quad t\ge 0,
\end{align*}
where $\mathfrak{A}_L:D(\mathfrak{A}_L)\subset \mathcal{X}\to \mathcal{X}$ is defined in \eqref{Generator}. This technique is originating
in \cite{Mil-75} has been widely investigated by many authors, in order to study the well-posedness of free-delay integro-differential equation. More references can be found in \cite{Pruss93}. According to this transformation, we say that the delay-integro-differential equation \eqref{problem} is well-posed if and only if the operator $\mathfrak{A}_L$ generates a strongly continuous semigroup $(\mathfrak{T}_L(t))_{t\ge 0}$ on $\mathcal{X}$. To prove that $\mathfrak{A}_L$ is a generator on $\mathcal{X}$, we first adopt the decomposition
\begin{align*}
\mathfrak{A}_L=\mathfrak{A}+\mathcal{L},
\end{align*}
where the operators $\mathfrak{A}$ and $\mathcal{L}$ are defined in \eqref{frakA} and \eqref{def-calL}, respectively. Second, we use a perturbation theorem in \cite{HMR} to prove that the operator $\mathfrak{A}$ generates a strongly continuous semigroup on $\mathcal{X}$ (see Theorem \ref{thm-generation-mathfrakA}). Third, based on regular linear systems we show that the operator $\mathcal{L}$ is a Miyadera-Voigt perturbation for $\mathfrak{A}$. Finally, we apply \cite[Chap.III]{EnNag} to deduce that $\mathfrak{A}_L$ is a generator on $\mathcal{X}$.

If we denote by $L_\Lambda$ the Yosida extension of $L$ with respect to the left shift semigroup on $L^p([-r,0],X)$ (see \eqref{Y-ext} for the definition), we prove that (see Theorem \ref{new-variation}) the solution of the above Cauchy problem (hence of the delay-integro-differential equation \eqref{problem}) is given by
\begin{align*}
z(t)=\begin{pmatrix}x(t)\\ v(t,\cdot)\\ x_t\end{pmatrix},\quad t\ge 0,
\end{align*}
where the function $t\mapsto x(t)$ is given by the following variation of parameters formula
\begin{align*}
x(t)=R(t)x+\int^t_0 R(t-s)(L_\Lambda x_s+f(s))ds,
\end{align*}
and $v(t,\cdot)$ is the solution of boundary system (\ref{input-shift}). Here $(R(t))_{t\ge 0}$ is the resolvent family associated the free-delay integro-differential equation problem (see. Definition (\ref{resolvent})).

Concerning the spectral theory for the generator  $\mathfrak{A}_L$, we have proved in Theorem \ref{thm-stability} that for $\lambda\in {\C}_0\cap\rho(\mathcal{A}^{\delta_0})\cap \rho(\mathcal{A}_0)$
\begin{align*}
\lambda\in\rho(\mathfrak{A}_L)\quad\Longleftrightarrow\quad \lambda\in\rho\Big(1+\hat{a}(\lambda))A+Le_\lambda\Big),
\end{align*}
where $\hat{a}$ is the Laplace transform of $a$ and $\mathcal{A}^{\delta_0}$ and $\mathcal{A}_{0}$ are the operators defined in \eqref{A-delta-0}. Remark that for $a(\cdot)\equiv0$ we obtain the same spectral result as in \cite{BatPia}. On the other hand, for $L\equiv 0,$ we identify the operators $\mathfrak{A}_L$ and $\mathcal{A}^{\delta_0},$ and then we retrieve the well-known results on the spectrum of $\mathcal{A}^{\delta_0}$, see e.g. \cite[Chap.7]{EnNag} and \cite{Pruss93}.

In the second part of this paper (see Section \ref{last-sect}), we prove that integrodifferential equations with state, input and output delays (see equation \eqref{problem-1}) form a subcalss of infinite-dimensional linear systems in the Salamon-Weiss sense \cite{WeiRegu}. This extends some existing results for standard delay and neutral systems \cite{Ha-Id-IMA}, \cite{HIR} and \cite{BH} respectively.

The organization of the paper is as follow: Section \ref{sec:2} is devoted to a brief background on feedback theory of regular linear systems in Salamon-Weiss sense. In Section \ref{sec:3} we prove the well-posedness of the problem \eqref{problem}. In Section \ref{sec:4} we propose a spectral study for state-delay integrodifferential equations.  The last section is concerned with reformulating integrodifferential equations as distributed systems.

\vspace {1cm}
\noindent{\bf Notation}\\
Throughout the paper we shall frequently use the following symbols:
Let $(X, ||\;||)$ be a Banach space and $\gamma$ be a real constant. For $p\in [1, \infty)$, we denote by
$L^p_\gamma(\R_+, X)$ stands for the space of functions $f: \R_+\to X$, such that
$e^{-\gamma t}f(t)$ is Bochner integrable on $\R_+$. $W^{1,p}(\R_+,X)$ is the Sobolev space associated to $L^p(\R_+, X):=L^p_0(\R_+, X)$.
Given two Banach spaces $X$ and $Y$, we denote the space of bounded
linear operators $X\to Y$ by $\L(X, Y)$, and for $G\in \L(X, Y),\; \|G\|_{\L(X,Y)}$ means
its operator norm (some time we only write $\|G\|$). We shall omit the subscripts if no confusion is possible. The identity operator will be denoted by $I$.
The Laplace transform of a function $f\in L^1_\omega(\R_+, X)$ is defined by $\hat{f}(\lambda)=\int_0^\infty
e^{-\lambda t}f(t)dt$ for $\lambda > \gamma$.  If $G:\R_+\to \L(X, Y)$ is an operator
valued function such that for each $x\in X$ the function $G(t)x$ is in $L^1_\omega(\R_+, X)$, then we define the linear operator $\hat{G}(\lambda)$ by $\hat{G}(\lambda)x = \int_0^\infty e^{-\lambda t}G(t)xdt$. Moreover, we shall use convolution of operator valued functions: Let $X,Y,Z$ be Banach spaces, and $G: \R_+\to\L(X, Y)$ and $F:\R_+\to L(Y,Z)$ functions such that $G(t)x$ and $F(t)x$ are measurable for all $x\in X$ or $Y$, respectively, and $||G(t)||, ||F(t)||$ are in $L^1_\gamma(\R_+, \R)$. Then $F\ast G: \R_+\to \L(X, Z)$ is defined by $(F\ast G)(t)x = \int_0^t F(t-s)G(s)xdt$ for $x\in X$. One can prove that $\widehat{(F\ast G)} = \hat{F}\cdot\hat{G}$. For a closed, linear operator
$A : D(A) \subset X \to X$, the resolvent set of $A$ is given by $\rho(A):=\{\lambda\in \C:\lambda-A;D(A) \subset X \to X\;\hbox{ is bijective}\}$.

Let $A:D(A)\subset X\to X$ be a generator of a strongly continuous semigroup $\T:=(T(t))_{t\ge 0}$ on $X$. We introduce a new norm $\|x\|_{-1}:=\|R(\alpha,A)x\|$ for $x\in X$ and some (hence all) $\alpha\in\rho(A)$. The completion of $X$ with respect to the norm $\|\cdot\|_{-1}$ is a Banach space denoted by $X_{-1}$. Moreover, the semigroup $\T$ can be extended to another strongly continuous semigroup $\T_{-1}:=(T_{-1}(t))_{t\ge 0}$ on $X$, whose generator $A_{-1}:X\to X_{-1}$ is the extension of $A$ to $X$. For more details on extrapolation theory we refer to \cite[ChapII]{EnNag}.

\section{Some background on regular linear systems} \label{sec:2}
In this section, we briefly recall the concept of regular linear systems (for more details on this theory we refer to \cite{WeiRegu}. Let $X,U$ and $Z$ be Banach spaces such that $Z\subset X$ densely and continuously. Let $A_m:Z\to X$ be a closed operator on $X,$ $G:Z\to U$ a linear operator (trace operator), and consider the boundary control problem
\begin{align}\label{BCP}
\begin{cases} \dot{z}(t)=A_m z(t),\quad z(0)=x,& t> 0,\cr  Gz(t)=u(t),& t\ge 0,\end{cases}
\end{align}
for $x\in X$ and $u\in L^p([0,+\infty),U)$. Under the following conditions
\begin{itemize}
  \item [{\bf(H1)}] $A:=A_m$ with domain $D(A)=\ker(G)$ generates a strongly continuous semigroup $T=(T(t))_{t\ge 0}$,
  \item  [{\bf(H2)}] $G$ is surjective,
\end{itemize}
the inverse $
D_\lambda:=\left( G_{|\ker(\lambda-A_m)}\right)^{-1}\in\mathcal{L}(U,X),\quad \lambda\in \rho(A),
$
exists, see \cite{Gr}. In addition if we seclect
$B:=(\lambda-A_{-1})D_\lambda,\;\la\in\rho(A)
$, then we have  $A_m=(A_{-1}+BG)_{|Z}$.
This implies that the boundary system \eqref{BCP} can be reformulated as a distributed system of the form
\begin{align}\label{DCS} \begin{cases} \dot{z}(t)=A_{-1}z(t)+Bu(t), & t\ge 0,\cr z(0)=z.\end{cases}
\end{align}
The mild solution of \eqref{DCS} is given by
\begin{align}\label{expression-mild}
\begin{split}
z(t)&=T(t)x+\int^{t}_0 T_{-1}(t-s)Bu(s)ds\cr &=T(t)x+\Phi^{A,B}_t u
\end{split}
\end{align}
for any $t\ge 0,$ $x\in X$ and $u\in L^p([0,+\infty),U)$. This solution takes its values in $X_{-1}$. It is then more practice to consider control operators $B\in \mathcal{L}(U,X_{-1})$ for which the solution $z(t)\in X$ for any $t\ge 0,$ initial condition $z\in X$ and control function $u\in L^p([0,+\infty),U)$. This situation happens if there exists $\tau>0$ such that $\Phi^{A,B}_{\tau} u\in X$ for any $u\in L^p([0,+\infty),U)$, see \cite{Weiss-control-89}. In fact, this condition also implies that $\Phi^{A,B}_t\in\mathcal{L}(L^p([0,+\infty),U),X)$ for any $t\ge 0$, due to the closed graph theorem. In this case, we say that $B$ is an admissible control operator for $A$.

 Let us now consider the observation function
\begin{align}\label{OS}
y(t)=\mathscr{C} z(t),& t\ge 0,
\end{align}
where $t\mapsto z(t)$ is the solution of \eqref{DCS}  and  $\mathscr{C}:Z\to U$ is a linear operator.

The system \eqref{DCS}-\eqref{OS} is called well-posed if $t\mapsto y(t)$ can be extended to a function $y\in L^p_{loc}([0,+\infty),U)$ that satisfies
\begin{align}\label{y-estimation}
\|y\|_{L^p([0,\al],U)}\le c \left(\|x\|+\|u\|_{L^p([0,\al],U)}\right)
\end{align}
for any $x\in X,\;u\in L^2_{loc}([0,+\infty),U),$ and some constants $\al>0$ and $c:=c(\al)>0$. Let us now introduce conditions on $B$ and $\mathscr{C}$ that guaranties the well-posedness of the system \eqref{DCS}-\eqref{OS}. To this end, first consider the operator
 \begin{align*}
 C:=\mathscr{C}_{|D(A)}\in\mathcal{L}(D(A),U).
\end{align*}
We say that $C$ is an {\em admissible observation} operator for $A$ if for some (hence all) $\alpha >0$, there exists $\kappa$ such that
\begin{align*}
\int^{\alpha}_0 \|C T(t)x\|^p dt\le \kappa^p \|x\|^p
\end{align*}
for all $x\in D(A).$ In this case, we have $T(t)x\in D(C_{\Lambda})$ for all $x\in X$ and a.e. $t>0$ and
\begin{align*}
\int^{\alpha}_0 \|C_{\Lambda} T(t)x\|^p dt\le \gamma^p \|x\|^p,\quad \forall x\in X,
\end{align*}
where $C_{\Lambda}$ is the Yosida extension of $C$ with respect to $A$ defined as follow
\begin{align}\label{Y-ext}\begin{split}
D(C_{\Lambda})&:=\{x\in X:\lim_{s\to +\infty} s C R(s,A)x\;\text{exists in}\; U\}\cr
C_{\Lambda} x&:=\lim_{s\to +\infty} s C R(s,A)x.
\end{split}
\end{align}
We need the following space
\begin{align*}
W^{2,p}_{0,t}(U):=\left\{u\in W^{2,p}([0,t],U):u(0)=0\right\},\qquad t>0.
\end{align*}
which is dense  in $L^p([0,t],U)$. Now a simple integration by parts yields $\Phi^{A,B}_t u\in Z$ for any $t\ge 0$ and $u\in W^{2,p}_{0,t}(U)$. This allows us to define  the so-called input-output map
\begin{align*}
(\mathbb{F} u)(t)=\mathscr{C}\Phi^{A,B}_t u,\qquad t\ge 0,\;u\in W^{2,p}_{0,t}(U).
\end{align*}
Using \eqref{expression-mild}, for  $x\in D(A)$ and $u\in W^{2,p}_{0,\tau}(U)$ ($\tau> 0$), the output function $y$ defined in \eqref{OS} satisfies
\begin{align*}
y(t)=CT(t)x+(\F u)(t)
\end{align*}
for any $t\in [0,\tau]$ and $u\in L^p([0,\tau],U)$.
\begin{defi}\label{ABC-well-posed}
Let $A,B,C$ and $\mathbb{F}$ as above. We say that the operator triple $(A,B,C)$  is well-posed if the following assertions hold:
\begin{itemize}
  \item [{\rm (i)}] $B$ is an admissible control operator for $A$,
  \item [{\rm (ii)}] $C$ is an admissible observation operator for $A$, and
  \item [{\rm (iii)}] There exist $\tau>0$ and $\kappa>0$ such that
  \begin{align*}
  \left\|\mathbb{F}u \right\|_{L^p([0,\tau],U)}\le \kappa \left\|u \right\|_{L^p([0,\tau],U)}
  \end{align*}
  for all $u\in W^{2,p}_{0,\tau}(U)$.
\end{itemize}
\end{defi}
If the triple $(A,B,C)$ is well-posed, by density of $ W^{2,p}_{0,\tau}(U)$ in $L^p([0,\tau],U)$, we can extend $\mathbb{F}$ to a bounded operator in $\mathcal{L}(L^p([0,\tau],U))$, for any $\tau>0$. In particular, the estimation of the observation \eqref{y-estimation} holds, so that the system \eqref{DCS}-\eqref{OS} is well-posed.
\begin{defi}\label{ABC-regular}
A well-posed triple $(A,B,C)$ is called regular (with feedthrough zero) if the following limit exists
\begin{align*}
\lim_{\tau \to 0} \frac{1}{\tau} \int^{\tau}_0 \left(\F (\1_{\R^+}\cdot v)\right)(\sigma)d\sigma=0,
\end{align*}
for any $v\in U$.
\end{defi}
\begin{remark}\label{bounded-C} \begin{itemize}
                                  \item[{\rm (i)}] According to \cite{WeiRegu}, if $(A,B,C)$ is a regular triple, then $\Phi^{A,B}_t u\in D(C_\Lambda)$ and $(\F u)(t)=C_\Lambda \Phi^{A,B}_t u$ for a.e. $t\ge 0$ and all $u\in L^p_{loc}([0,+\infty),U)$. This shows also that the state trajectory and the output function of the system \eqref{DCS}-\eqref{OS} satisfy $x(t)\in D(C_\Lambda)$ and $y(t)=C_\Lambda z(t)$ for a.e. $t>0$, and all initial state  $z(0)=x\in X$ and all input $u\in L^p_{loc}([0,+\infty),U)$.
                                  \item [{\rm (ii)}] Clearly, if one of the operators $B$ or $C$ is bounded and the other is admissible for $A$, then the triple $(A,B,C)$ is regular.
                                \end{itemize}
\end{remark}
We end this section with an example of regular linear system which will be frequently used in some proofs in the next sections.
\begin{example}\label{regular-shift-system}
Let $X$ be a Banach space and $p>1$ be a real number.  It is known (see e.g. \cite[chap.2]{EnNag}) that the family $(S(t))_{t\ge 0}$ defined by
\begin{align}\label{shift-sg}
(S(t)\varphi)(\theta):=\begin{cases} 0,& t+\theta\ge 0,\cr \varphi(t+\theta),& t+\theta\le 0\end{cases}
\end{align}
for any $f\in L^p([-r,0],X)$, $t\ge 0$ and $\theta\in [-r,0]$,  is a strongly continuous semigroup on $L^p([-r,0],X),$ (called the left shift semigroup). The generator of this semigroup is
\begin{align}\label{generator-left-shift}
Q\varphi=\varphi',\quad D(Q)=\left\{\varphi\in W^{1,p}([-r,0],X):\varphi(0)=0\right\}.
\end{align}
 Now, consider the boundary system
\begin{align}\label{input-shift}
\begin{cases}
\displaystyle{\frac{\partial v(t,\theta)}{\partial t}}=\frac{\partial v(t,\theta)}{\partial \theta},& t\ge 0,\; \theta\in [-r,0],\cr v(0,\theta)=\varphi(\theta),& \theta\in [-r,0],\cr v(t,0)=x(t),& t\ge 0,
\end{cases}
\end{align}
We select $Q_m:=\displaystyle{}\frac{\partial}{\partial \theta}$ with maximal domain $D(Q_m)=W^{1,p}([-r,0],X)$, so $Q=Q_m$ and $D(Q)=\ker G$ with $G=\delta_0$ where $\delta_0f=f(0)$ is the Dirac operator. The Dirichlet operator $d_{\lambda}$ associated to \eqref{input-shift} is given by
\begin{align*}
d_\lambda x=e_\lambda x:=e^{\lambda \cdot}x,\quad \lambda\in\rho(Q)=\mathbb{C},\;x\in X.
\end{align*}
We put
\begin{align*}
\beta:=(\lambda-Q_{-1})d_\lambda,\qquad \lambda\in\mathbb{C}.
\end{align*}
The control maps associated with the control operator $\beta$ are given by
\begin{align*}
\left(\Phi^X_t u\right)(\theta)=\begin{cases} u(t+\theta),& -t\le \theta\le 0,\cr 0,& -r\le \theta<-t,\end{cases}
\end{align*}
for any $t\ge 0$ and $u\in L^p([0,+\infty),X)$, see \cite{HIR}. Thus the operator $\beta\in\mathcal{L}(X,(L^p([-r,0],X))_{-1})$ is an admissible control operator for $Q$. Now consider the operator $L:W^{1,p}([-r,0],X)\to X$ defined by \eqref{L-RS} and  define
\begin{align*}
L_0:=L_{|D(Q)}\in \mathcal{L}(D(Q),X).
\end{align*}
Then $L_0$ is an admissible observation operator for $Q$, see \cite[Lemma 6.2.]{Hadd-SF}. We select
\begin{align*}
(\F^X u)(t)=L \Phi_t u,\quad u\in W^{2,p}_{0,t}(X).
\end{align*}
As in \cite{HIR}, we show that for any $\tau>0$ and $u\in W^{2,p}_{0,\tau}(X)$,
\begin{align*}
\int^{\tau}_0 \|(\F^X u)(t)\|^pdt & \le \int^{\tau}_0 \left(\int^0_{-t}\|u(t+\theta)\| d|\mu|(\theta)\right)^p dt\cr &\le \left(|\mu|([-\tau,0])\right)^{\frac{p}{q}}\int^{\tau}_0 \int^0_{-t}\|u(t+\theta)\|^p d|\mu|(\theta)dt
\cr &\le \left(|\mu|([-\tau,0])\right)^{\frac{p}{q}}\int^0_{-\tau} \int^{\t-\theta}_0\|u(t+\theta)\|^p dtd|\mu|(\theta)\cr & \le \left(|\mu|([-\tau,0])\right)^{p}\int^\tau_0 \|u(\si)\|^pd\si,
\end{align*}
due to H\"older inequality and Fibini's theorem. This shows that the triple $(Q,\beta,L_0)$ is well-posed. On the other hand, it is shown in \cite[Theorem 3]{HIR} that the triple $(Q,\beta,L_0)$ is a regular linear system.
\end{example}

\section{Well-posedness of the state-delay integro-differential equation}\label{sec:3}
In this section, we will study the well-posedness of the state-delay integro-differential equation \eqref{problem}. But first we recall some facts about the following free-delay integro-differential equation
\begin{align}\label{free-voltera}
\dot{x}(t)=Ax(t)+\int^t_0 a(t-s)Ax(s)ds+f(t),\quad x(0)=x,\quad t>0.
\end{align}
It is known that the solvability of \eqref{free-voltera} is related to the concept of strongly continuous
family defined  as follows:
\begin{defi}\label{resolvent}
Let $a\in L_{loc}^{1}\left(\mathbb{R}^{+}\right)$. A strongly continuous
family $\left(R\left(t\right)\right)_{t\geq0}\subset\mathcal{L}\left(X\right)$ is called resolvent
family for the homogeneous free-delay integro-differential equation, if the following three
conditions are satisfied:
\begin{itemize}
\item $R(0)=I$.
\item $R(t)$ commutes with $A,$ which means $R(t)D(A)\subset D(A)$
for all $t\geq0,$ and $AR(t)x=R(t)Ax$ for all $x\in D(A)$ and $t\geq0$.
\item For each $x\in D(A)$ and all $t\geq0$ the resolvent equations
holds:
\end{itemize}
\[
R(t)x=x+\int_{0}^{t}a\left(t-s\right)AR(s)xds.
\]

\end{defi}
Generally, in the definition above, it is not necessary that A be the generator of a semigroup. It is well-known (see, e.g., [7, Section 1]) that the homogeneous free-delay integro-differential system is well-posed if and only if it has a resolvent $R(t)$. In this situation, $x(t)=R(t)x, t\ge 0$  with $x\in X$ is the mild solution, which gives the unique classical solution if $x\in D(A)$. For more details on resolvent families we refer to the monograph by Pruss \cite{Pruss93}.

It is well-known that the assumptions $a(\cdot)\in W^{1,p}(\R_+,\C)$ and $A$ generates a semigroup on $X$ imply that the free-delay integro-differential equation has a unique resolvent family $(R(t))_{t\ge 0}$, see e.g. \cite{DS-85}, \cite{Pruss93}.

In the sequel, we will use matrices operators to solve the equation \eqref{free-voltera}. We then
select
\begin{align*}
\mathcal{X}_0:=X\times L^p(\R^+,X)\quad\text{with norm}\quad \left\|(\begin{smallmatrix} x\\ f\end{smallmatrix})\right\|:=\|x\|+\|f\|_p.
\end{align*}
In this space, we consider the following unbounded operator matrices
\begin{align}\label{A-delta-0}
\begin{split}
&\mathcal{A}_{0}:=\begin{pmatrix}A&\delta_0\\0&\displaystyle{\frac{d}{ds}} \end{pmatrix},\quad \mathcal{A}^{\delta_0}:=\begin{pmatrix}A&\delta_0\\a(\cdot)A&\displaystyle{\frac{d}{ds}} \end{pmatrix}\cr
& D(\mathcal{A}^{\delta_0})=D(\mathcal{A}_{0}):=D(A)\times W^{1,p}(\mathbb{R}_+,X),
\end{split}
\end{align}
where  $\displaystyle{\frac{d}{ds}}$ is the first derivative with $D(\displaystyle{\frac{d}{ds}})=W^{1,p}(\mathbb{R}^{+},X)$.

The following result is well known, see e.g. \cite{EnNag}, \cite [p.339]{Pruss93} and \cite{FadiliBounitb}.
\begin{theorem}
  The operator $\mathcal{A}^{\delta_0}$ generates a strongly continuous semigroup $(T^{\delta_0}(t))_{t\geq 0}$ on $\mathcal{X}_0$. Moreover, if we assume that $a(\cdot)\in W^{1,p}(\R_+,\C)$, then \begin{align*}
T^{\delta_0}(t)=\begin{pmatrix}
R(t)&\Upsilon(t)\\ \ast&\ast
\end{pmatrix},\quad t\ge 0,
\end{align*}
with
\begin{align}\label{resolvent-Volterra}
\Upsilon(t)f:=\int_0^tR(t-s)f(s)ds,\quad f\in L^p(\R^+,X).
\end{align}
In particular $R(t)$ is exponentially bounded since $(T^{\delta_0}(t))_{t\ge 0}$ is so.
\end{theorem}
Let us now solve the state-delay integro-differential equation \eqref{problem}. The latter can be reformulated in $\mathcal{X}_0$ as the following problem
\begin{align} \label{auxillary-problem}
\begin{cases}
\dot{\varrho}(t)=\mathcal{A}^{\delta_0}\varrho(t)+\Big(\begin{smallmatrix} Lx_t\\ 0\end{smallmatrix}\Big),& t\ge 0,\cr
\varrho(0)=(\begin{smallmatrix} x\\ f\end{smallmatrix}),\quad x_0=\varphi.
\end{cases}
\end{align}
This equation is not yet a Cauchy problem because we still have a delay term. In order to reformulate it as Cauchy problem, we need the following larger product state space
\begin{equation}\label{space-calX}
\mathcal{X}:=\mathcal{X}_0\times L^p([-r,0],X),
\end{equation}
and the new state
\begin{align}\label{function-z}
z(t):=\begin{pmatrix}
\varrho(t)\\ x(t+\cdot)
\end{pmatrix},\qquad t\ge 0.
\end{align}
By combining \eqref{input-shift} and  \eqref{auxillary-problem}, we rewrite the problem \eqref{problem}  as the following boundary value problem
\begin{align}\label{bound}
\begin{cases}
\dot{z}(t)=\mathcal A_{m,L}z(t),&t\geq 0\cr
z(0)=z^0,\\
Gz(t)=Mz(t),& t\ge 0,
\end{cases}
\end{align}
where the operator $\mathcal A_{m,L}:\mathcal{Z}\to \mathcal{X}$ is given by
\begin{align*}
\mathcal A_{m,L}:= \left(
\begin{array}{c|c}
 \mathcal{A}^{\delta_0}&\begin{matrix}
L\\0
\end{matrix}\\
\hline
\begin{matrix}
0&0
\end{matrix}& \frac{d}{d\theta}
\end{array}\right),\quad \mathcal{Z}:=D(\mathcal{A}^{\delta_0})\times W^{1,p}([-r,0],X),
\end{align*}
and $G:\mathcal{Z}\to X$ and $M:\mathcal{X}\to X$ are the linear operators
\begin{align*}
G:=\begin{bmatrix}
0&0&\delta_0
\end{bmatrix},\quad M:=\begin{bmatrix}
I&0&0
\end{bmatrix},
\end{align*}
and the initial state
\begin{align*}
z^0=\left(\begin{smallmatrix}x\\f\\ \varphi\end{smallmatrix}\right).
\end{align*}
To solve the problem \eqref{problem}, it suffices to solve the following Cauchy problem
\begin{align}\label{perturbed-equation}
\begin{cases}
\dot{z}(t)=\mathfrak{A}_L z(t),& t\ge 0,\cr z(0)=z^0,
\end{cases}
\end{align}
where
\begin{align}\label{Generator}
\mathfrak{A}_L:=\mathcal{A}_{m,L},\qquad D(\mathfrak{A}_L):=\left\{z\in \mathcal{Z}: Gz=Mz\right\}.
\end{align}
This means that it suffices to show that the operator $\mathfrak{A}_L$  is a generator of a strongly continuous semigroup on $\mathcal{X}$. Thus, we define
\begin{align}\label{def-calL}
\mathcal A_{m,0}:= \left(
\begin{array}{c|c}
 \mathcal{A}^{\delta_0}&\begin{matrix}
0\\0
\end{matrix}\\
\hline
\begin{matrix}
0&0
\end{matrix}& \displaystyle\frac{d}{d\theta}
\end{array}\right) \quad\text{and}\quad \mathcal{L}:=\left(
\begin{array}{c|c}
0&\begin{matrix}
L\\0
\end{matrix}\\
\hline
\begin{matrix}
0&0
\end{matrix}& 0
\end{array}\right)
\end{align}
on $D(\mathcal A_{m,L})=\mathcal{Z}$, and
\begin{align}\label{frakA}
\mathfrak{A}:=&\mathcal{A}_{m,0},\cr D(\mathfrak{A}):=&\left\lbrace z\in \mathcal{Z},\quad Gz=Mz\right\rbrace.
\end{align}
Observe that we have
\begin{align*}
\mathfrak{A}_L=\mathfrak{A}+\mathcal{L}.
\end{align*}
Next we follow the following strategy: We first prove that $\mathfrak{A}$ is a generator on $\mathcal{X}$ and second show that $\mathcal{L}$ is a Miyadera-Voigt perturbation for $\mathfrak{A}$.

Let us start by proving that $\mathfrak{A}$ is a generator. In fact, define the operator
\begin{align*}
\mathcal A := \mathcal{A}_{m,0}, \text{ with domain }D(\A)=\ker(G).
\end{align*}
It is clear that the operator $\mathcal{A}$ generates a strongly continuous semigroup $(\T(t))_{t\geqslant 0}$ on $\mathcal{X}$ given by
\begin{equation}\label{Volterra}\T(t)\begin{pmatrix}
x\\f\\\varphi
\end{pmatrix}=\begin{pmatrix}
T^{\delta_0}(t)\begin{pmatrix}
x\\f
\end{pmatrix}\\
\hline
S(t)\varphi
\end{pmatrix},\quad t\geq 0,\;\begin{pmatrix}
x\\f\\\varphi
\end{pmatrix}\in\mathcal{X}.\end{equation}
Now let us compute $\mathbb D_\lambda$, the Dirichlet operator associated with $G$ and $\mathcal{A}_{m,0}$. To determine $\mathbb D_\lambda$, it suffices to determine
Ker$(\lambda-\mathcal{A}_{m,0})$, for $\lambda\in\rho(\mathcal{A}^{\delta_0})$.
For $\lambda\in\rho(\mathcal{A})=\rho(\mathcal{A}^{\delta_0})$,
\begin{eqnarray*}
\hbox{Ker}(\lambda-\mathcal{A}_{m,0})&=\hbox{Ker}(\lambda-\mathcal{A}^{\delta_0})\times \hbox{Ker}(\lambda-Q_m)\cr
&=\left\{\Big(\begin{smallmatrix}0\\0\end{smallmatrix}\Big)\right\}\times \Big\{e_\lambda v: v\in X\Big\}.
\end{eqnarray*}
This implies that
\begin{align*}
\mathbb D_\lambda v=\begin{pmatrix}0\\0
\\ e_\lambda v
\end{pmatrix}, \quad \lambda\in \rho(\mathcal{A}^{\delta_0}),\; v\in X.
\end{align*}
Now, let us consider the operator
\begin{align}\label{calB}
\mathcal{B}:=(\lambda-\mathcal{A}_{-1})\mathbb{D}_\lambda\in\mathcal{L}(X,\mathcal{X}_{-1}),\quad \lambda\in \rho(\mathcal{A}).
\end{align}

The facts mentioned above give the following variation of the formula of the parameters which will be useful later.

\begin{theorem}\label{thm-generation-mathfrakA}
The operator $\mathfrak{A}$ generates a strongly continuous semigroup $(\mathfrak{T}(t))_{t\ge 0}$ on $\mathcal{X}$ satisfying
\begin{equation}\label{FVC-mathfrak-T}
\mathfrak{T}(t)z^0=\T(t)z^0+\int_0^t\T_{-1}(t-s)\mathcal{B}M\mathfrak{T}(s)z^0 ds
\end{equation}
for any $z^0\in\mathcal{X}$ and $t\ge 0$.
\end{theorem}
\begin{proof}
For $u\in L^p([0,\infty),X)$ and $\lambda>0$ sufficiently large,  using the Laplace-transform we obtain
\begin{eqnarray*}
\widehat{(\mathbb{T}_{-1}\ast\mathcal{B} u)}(\lambda)=R(\lambda,\mathcal{A}_{-1})\mathcal{B}\widehat{u}(\lambda)
=
\mathbb{D}_\lambda \hat u(\lambda)=\begin{pmatrix}0\\0
\\ e_\lambda \hat u(\lambda)
\end{pmatrix}.
\end{eqnarray*}
Thus
\begin{eqnarray*}
\widehat{(\mathbb{T}_{-1}\ast\mathcal{B}u} )(\lambda)={\begin{pmatrix}0\\0
\\ \widehat{\Phi_\cdot^{Q,\beta} u}(\lambda)
\end{pmatrix}},
\end{eqnarray*}
where $(\Phi_t^{Q,\beta})_{t\ge 0}$ is the control maps associated to the regular system $(Q,\beta,L_0)$ defined before in Example \ref{regular-shift-system}.
By injectivity of Laplace-transform, we deduce that
\begin{equation}
\int_0^t\T_{-1}(t-s)\mathcal{B}u(s)ds=\begin{pmatrix}
0\\0\\
\Phi_t^{Q,\beta} u
\end{pmatrix}\in \mathcal{X}.
\end{equation}
It follows that $\mathcal{B}$ is an admissible control operator for $\mathcal{A}$. The fact that $M\in \mathcal L(\mathcal{X},X)$,
$(\mathcal{A},\mathcal{B},M)$ is  a regular system with $I:X\to X$ as an admissible feedback (see. \cite{WeiRegu}). Therefore, the operator
\begin{align}
\mathcal{A}^{cl}:=&\mathcal{A}_{-1}+\mathcal{B}M\cr
D(\mathcal{A}^{cl}):=&\left\lbrace z\in\mathcal{X},\; \mathcal{A}^{cl}z\in\mathcal{X}\right\rbrace
\end{align}
generates a strongly continuous semigroup $(T^{cl}(t))_{t\ge 0}$ on $\mathcal{X}$ satisfying
\begin{equation*}
T^{cl}(t)z^0=\T(t)z^0+\int_0^t\T_{-1}(t-s)\mathcal{B}MT^{cl}(s)z^0 ds
\end{equation*}
for any $z^0\in \mathcal{X}$ and $t\ge 0$, see \cite{WeiRegu}. Finally, we mention that the operator $\mathcal A^{cl}$ coincides with the operator $\mathfrak{A}$, due to \cite[Theorem 4.1]{HMR}.
\end{proof}
The following theorem is the main results of this section.

\begin{theorem}\label{generation-big-A}
The operator \begin{align*}
\mathfrak{A}_{L}:=\mathfrak{A}+\L,\quad
D(\mathfrak{A}_{L}):=D(\mathfrak{A})
\end{align*}
 generates a strongly continuous semigroup $(\mathfrak{T}_{L}(t))_{t\geq 0}$ on $\mathcal{X}$.
\end{theorem}
\begin{proof}
To prove that $\mathfrak{A}_{L}$ is a generator it suffices to show that $\L$ is an admissible observation operator for $\mathfrak{A}$ (see. \cite[Theorem 2.1]{Hadd-SF}). To this end, let us define the operator
\begin{align*}
\L_{0}=\L_{|D(\mathcal{A})}\in \mathcal{L}(D(\mathcal{A}),\mathcal{X}).
\end{align*}
Let us first show that $\L_0$ is an admissible observation
operator for $\A$. For $\alpha>0$ and $z^0=(x,f,\varphi)^\top\in D(\A)=D(\mathcal{A}^{\delta_0})\times D(Q),$ we have
\begin{align*}
\int^\alpha_0 \|\L_0 \mathbb{T}(t)z^0\|^pdt=\int^\alpha_0 \|L_0 S(t)\varphi\|^pdt,
\end{align*}
where $L_0$ is the operator defined in Example \ref{regular-shift-system}, which is an admissible observation operator for $Q$. It follows from the above
estimate that $\L_0$ is admissible for $\A$. On the other hand, by using a similar argument as in \cite[Lemma 6.3]{Hadd-SF},
one can see that the Yosida extension of $\L_0$ with respect to $\mathcal{A}$ is explicitly given by
\begin{align}\label{Lgamma}
D(\L_{0,\Lambda})=\mathcal{X}_0\times D(L_{0,\Lambda}),\qquad \L_{0,\Lambda}=\left(\begin{array}{c|c}
0&\begin{matrix}
L_{0,\Lambda}\\0
\end{matrix}\\
\hline
\begin{matrix}
0&0
\end{matrix}& 0
\end{array}\right).
\end{align}
In addition for any $z^0\in D(\mathfrak{A})$, we have
\begin{align}\label{calL-0-estimcalL-0-Ga}
\int^\alpha_0 \|\L_{0,\Lambda} \mathbb{T}(t)z^0\|^pdt\le c^p \|z^0\|^p
\end{align}
for a constant $c>0$. According to \cite[Lemma 3.6]{HMR}, we have
\begin{align*}
\mathcal{Z}\subset D(\L_{0,\Lambda}),\quad\text{and}\quad \left(\L_{0,\Lambda}\right)_{|\mathcal{Z}}=\L.
\end{align*}
Thus for $z^0\in D(\mathfrak{A})$,
\begin{align}\label{egalite-L-L}
\int^\alpha_0 \|\L \mathfrak{T}(t)z^0\|^pdt=\int^\alpha_0 \|\L_{0,\Lambda} \mathfrak{T}(t)z^0\|^pdt.
\end{align}
If we consider
\begin{align*}
u(t):=M\mathfrak{T}(t)z^0,\qquad t\ge 0,
\end{align*}
then by using the proof of Theorem \ref{thm-generation-mathfrakA}, we obtain
\begin{align}\label{varcons-T}
\int_0^t\T_{-1}(t-s)\mathcal{B}MT^{cl}(s)z^0 ds=\bpm 0\\0\\\Phi^{Q,\beta}_{t}u\epm.
\end{align}
Now from Example \ref{regular-shift-system}, $(Q,\beta,L_0)$ is regular and that $\Phi^{Q,\beta}_{t}u\in D(L_{0,\Lambda})$ for a.e. $t\ge 0$. Combining (\ref{Lgamma}) and (\ref{varcons-T}) we obtain
\begin{align*}
\int_0^t\T_{-1}(t-s)\mathcal{B}MT^{cl}(s)z^0 ds\in D(\L_{0,\Lambda}),
\end{align*}
and
\begin{align}\label{F-estim}
\begin{split}
\int^\alpha_0 \left\|\L_{0,\Lambda}\int_0^t\T_{-1}(t-s)\mathcal{B}MT^{cl}(s)z^0 ds\right\|^pdt&
=\int^\alpha_0\|L_{0,\Lambda}\Phi^{Q,\beta}_tu\|^pdt\cr & \le \kappa^p \|u\|^p_{L^p([0,\alpha],X)}\cr &\le \tilde{\kappa}^p \|z^0\|^p,
\end{split}
\end{align}
for a constant $\tilde{\kappa}>0;$ due to the facts $M$ is bounded and that $\mathfrak{T}$ is exponentially bounded.
Finally, using on the one side, \eqref{FVC-mathfrak-T} and \eqref{egalite-L-L}, on the other side, the estimates \eqref{calL-0-estimcalL-0-Ga} and
\eqref{F-estim}, we obtain
\begin{align*}\label{egalite-L-L}
\int^\alpha_0 \|\L \mathfrak{T}(t)z^0\|^pdt\le \gamma^p \|z^0\|^pdt,
\end{align*}
for any $z^0\in D(\mathfrak{A})$ and some constant $\gamma:=\gamma(\alpha)>0$. This ends the proof.
\end{proof}

\begin{remark} \label{spactrum}
Observe that $M\D_\la=0$ for any $\la\in \rho(\mathcal{A})$. It is shown in \cite[Theorem 4.1]{HMR} that
for $\lambda\in \rho(\A)$ we have
$$\lambda\in\rho(\mathfrak{A})\quad\Longleftrightarrow\quad 1\in\rho(\D_\lambda M)\quad\Longleftrightarrow\quad 1\in\rho(M\D_\lambda )=\C^*.$$
Thus
\begin{equation}
\rho(\A)=\rho(\mathcal{A}^{\delta_0})\subset \rho(\mathfrak{A}).
\end{equation}
Again by \cite[Theorem 4.1]{HMR}, for $\lambda\in\rho(\A)$, we have
$$R(\lambda,\mathfrak{A})=(I-\D_\lambda M)^{-1}R(\lambda,\A).$$
\end{remark}

In the following result, we give a new expression of the solution of the Cauchy problem \eqref{perturbed-equation} by appealing only to semigroup $(\mathbb{T}(t))_{t\ge 0}$.
\begin{prop}\label{sol}
The solution $z(\cdot)$ of the Cauchy problem \eqref{perturbed-equation} satisfies
\begin{align}\label{hamid}
\begin{split}
& z(s)\in D(\L_{\Lambda}),\quad a.e.\; s\ge 0,\cr
&z(t)=\T(t)z(0)+\int_0^t\T_{-1}(t-s)\B Mz(s)ds+\int_0^t\T(t-s)\L_{\Lambda}z(s)ds
\end{split}
\end{align}
 for any $t\ge 0,\;z(0)\in X,$
where $\mathcal{L}_\Lambda$ is the Yosida extension of $\mathcal{L}$ with respect to $\mathfrak{A}$.
\end{prop}
\begin{proof}Using \cite[Theorem 5.1]{Hadd-SF} and the fact that $\mathcal{L}$ is an admissible observation operator for $\mathfrak{A}$
(see the proof of the previous result), the solution of the Cauchy problem \eqref{perturbed-equation} satisfies
$z(s)\in D(\L_{\Lambda})$ for a.e.\; $s\ge 0$, and
$$z(t)=\mathfrak{T}(t)z(0)+\int^t_0\mathfrak{T}(t-s)\L_\Lambda z(s) ds,$$ for any $t\ge 0,\;z(0)\in X$.\\
 For a sufficiently large $\lambda>0,$ the Laplace transform of $z$ yields
\begin{align*}\widehat z(\lambda)&= R(\lambda,\mathfrak{A})z(0)+R(\lambda,\mathfrak{A})\L_\Lambda\widehat z(\lambda)]
\cr &= (I-\mathbb{D}_\lambda M)^{-1}\left[R(\lambda,\A)z(0)+R(\lambda,\A)\L_\Lambda \widehat z(\lambda)\right].\end{align*}
This in turn implies that
\begin{align}\label{E1}\widehat z(\lambda)=R(\lambda,\A)z(0)+\mathbb{D}_\lambda M \widehat z(\lambda)+R(\lambda,\A)\L_\Lambda \widehat z(\lambda).
\end{align}
Let
\begin{align*}
\varpi(t):=\T(t)z(0)+\int_0^t\T_{-1}(t-s)\B Mz(s)ds+\int_0^t\T(t-s)\L_{\Lambda}z(s)ds,\quad t\ge 0.
\end{align*}
According to \eqref{calB} and \eqref{E1}, the Laplace transform of $\varpi$ satisfies
\begin{align*}
\widehat \varpi(\lambda)&=
R(\lambda,\A)z(0)+R(\lambda,\A_{-1})\B M\widehat z(\lambda)+R(\lambda,\A)\L_\Lambda\widehat z(\lambda)\cr & =\widehat z(\lambda),\end{align*}
and the result follows by virtue of the injectivity the Laplace transform.
\end{proof}




Note that for $a(\cdot)\equiv0$,  the resolvent familly correspond to the semigroup generated by $A$. Therefore the following result which establishes the relation between the solutions of (\ref{problem}) and (\ref{perturbed-equation}), extends \cite[Proposition 2]{HIR} from delay differential equations to delay integro-differential one.

\begin{theorem}\label{new-variation}
For any initial condition $(x,f,\varphi)^\top\in\mathcal{X}$, there exists a unique  solution $x(\cdot)$ of (\ref{problem}) satisfying $ x_t\in D(L_\Lambda)$ for a.e. $t\ge 0$ and
\begin{align*}
x(t)=R(t)x+\int_0^tR(t-s)(L_\Lambda x_s+f(s)) ds,\quad t\ge 0.
\end{align*}
\end{theorem}

\begin{proof}
Let	$z(0)=(x,f,\varphi)^\top\in\mathcal{X}$ and $z(t)=(x(t),v(,\cdot),w(t))^\top$ be
the solution of the Cauchy problem \eqref{perturbed-equation} associated to $z(0)$. Appealing to Proposition (\ref{sol}) and by combining \eqref{hamid} together with \eqref{Volterra}, \eqref{varcons-T}, \eqref{resolvent-Volterra} we deduce that $w(t)= x_t$ and \eqref{perturbed-equation} becomes
\begin{align*}
\begin{pmatrix}
x(t)\\ v(t,\cdot)\\x_t
\end{pmatrix}&=\begin{pmatrix}
T^{\delta_0}(t)\begin{pmatrix}
x\\f
\end{pmatrix}\\
\hline
S(t)\varphi
\end{pmatrix}+\bpm 0\\0\\\Phi^{Q,\beta}x\epm+\int_0^t\begin{pmatrix}
T^{\delta_0}(t-s)\begin{pmatrix}
L_\Lambda x_s\\0
\end{pmatrix}\\
\hline
0
\end{pmatrix}ds\cr
&=
\begin{pmatrix}
\begin{pmatrix}
R(t)&\Upsilon(t)\cr \ast&\ast
\end{pmatrix}
\begin{pmatrix}
x\\ f
\end{pmatrix}
\\
\hline
S(t)\varphi
\end{pmatrix}
+\bpm 0\\0\\\Phi^{Q,\beta}x\epm+\int_0^t\begin{pmatrix}
\begin{pmatrix}
R(t-s)&\Upsilon(t-s)\cr \ast&\ast
\end{pmatrix}\begin{pmatrix}
L_\Lambda x_s\\0
\end{pmatrix}\\
\hline
0
\end{pmatrix}ds.
\end{align*}
Thus by equalizing the two first components of the above equality we get the required result.
\end{proof}
\section{Spectral theory}\label{sec:4}
In this section, we attempt to study the spectral theory for the generator $\mathfrak {A}_L $ and we get some results extending those obtained in \cite{BatPia} and \cite{EnNag}. First of all, we require the following straightforward lemma.
\begin{lem}\label{lemma-stability}
For $\lambda\in\rho(\mathcal{A}^{\delta_0})$, we have
\begin{equation}\label{res1}\lambda\in\rho(\mathfrak{A}_L)\quad\Longleftrightarrow\quad 1\in\rho(\L R(\lambda,\mathfrak{A})).
\end{equation}
\end{lem}
\begin{proof}
 For $\lambda\in\rho(\mathcal{A}^{\delta_0})$, we write $\lambda-\mathfrak{A}_L=\lambda-\mathfrak{A}- \L$. By Remark (\ref{spactrum}), we have $\lambda\in\rho(\mathfrak{A})$
\begin{eqnarray*}
\lambda-\mathfrak{A}_L=(1-\L R(\lambda,\mathfrak{A}))(\lambda-\mathfrak{A}).
\end{eqnarray*}
This implies that $\lambda\in\rho(\mathfrak{A}_L)$ if and only if $(1-\L R(\lambda,\mathfrak{A}))^{-1}$ exists.
\end{proof}

From the previous results we immediately obtain the following result which characterises the spectrum of $\mathfrak{A}_L$ and then extends \cite[Lemma 4.1]{BatPia} from delay differential or integro-differential equations to delay integro-differential equations (\ref{problem}) .

\begin{theorem}\label{thm-stability}
For $\lambda\in {\C}_0\cap\rho(\mathcal{A}^{\delta_0})\cap \rho(\mathcal{A}_0)$, we have
$$\lambda\in\rho(\mathfrak{A}_L)\quad\Longleftrightarrow\quad \lambda\in\rho\Big((1+\hat{a}(\lambda))A+Le_\lambda\Big).$$
\end{theorem}

\begin{proof}
By virtue of Remark (\ref{spactrum}), we have $\lambda\in\rho(\mathfrak{A})$ and therefore
\begin{eqnarray*}
R(\lambda,\mathfrak{A})&=&(I-\D_\lambda M)^{-1}R(\lambda,\A)\cr
&=&\bpm 1&0&0\\0&1&0\\ e_\lambda&0&1\epm\left(\begin{array}{c|c} R(\lambda,\mathcal{A}^{\delta_0})&\begin{matrix}0\\0\end{matrix}\\\hline \begin{matrix}0&0\end{matrix}&R(\lambda,Q)\\
\end{array}\right)\cr
&=&\left(\begin{array}{c|c} R(\lambda,\mathcal{A}^{\delta_0})&0\\\hline \bpm e_\lambda&0\epm R(\lambda,\mathcal{A}^{\delta_0})&R(\lambda,Q)\\
\end{array}\right).
\
\end{eqnarray*}
Appealing to \cite[Proposition. 7.25]{EnNag}, (see. \cite{FadiliBounitb}) leads to $\displaystyle{\frac{\lambda}{1+\hat{a}(\lambda)}}\in \rho(A)$ and by virtue of \cite[Lemma 7 (ii)]{FadiliBounitb}, we obtain
$$R(\lambda,\mathcal{A}^{\delta_0})
=\bpm H(\lambda)&
H(\lambda)\delta_0R(\lambda,\displaystyle{\frac{d}{ds}})\cr R(\lambda,\displaystyle{\frac{d}{ds}})
VH(\lambda)&R(\lambda,\displaystyle{\frac{d}{ds}})VH(\lambda)\delta_0R(\lambda,\displaystyle{\frac{d}{ds}})+R(\lambda,\displaystyle{\frac{d}{ds}})\epm,$$
where \\
$H(\lambda)= (\lambda-(1+\hat{a}(\lambda)A))^{-1}$ and $V\in \mathcal L(D(A), W^{1,p}(\R_+,X))$ is given by
$(Vx)(t):= a(t)Ax$ for $x\in D(A)$ and a.e. $t>0$.\\
Therefore, direct computations yield
$$R(\lambda,\mathfrak{A})=\bpm H(\lambda)& H(\lambda)\delta_0R(\lambda,\displaystyle{\frac{d}{ds}})&0\cr R(\lambda,\displaystyle{\frac{d}{ds}})VH(\lambda)&R(\lambda,\displaystyle{\displaystyle{\frac{d}{ds}}})VH(\lambda)\delta_0R(\lambda,\displaystyle{\frac{d}{ds}})+R(\lambda,\displaystyle{\frac{d}{ds}})&0\cr
e_\lambda H(\lambda)& e_\lambda H(\lambda)\delta_0R(\lambda,\displaystyle{\frac{d}{ds}})&R(\lambda,Q)\epm$$
and
\begin{equation}\label{LR}\L R(\lambda,\mathfrak{A})=\bpm Le_\lambda H(\lambda)& Le_\lambda H(\lambda)\delta_0R(\lambda,\displaystyle{\frac{d}{ds}})&LR(\lambda,Q)\cr
0&0&0\cr0&0&0\epm
\end{equation}
By combining \eqref{res1} and \eqref{LR}, we deduce that
$\lambda\in\rho(\mathfrak{A}_L)$ if and only if $1\in\rho(Le_\lambda H(\lambda))$. Rewriting
$$(\lambda-(1+\hat{a}(\lambda))A)-Le_\lambda=\left(1-Le_\lambda H(\lambda)\right)(\lambda-(1+\hat{a}(\lambda))A)$$
we conclude that $$\lambda\in\rho(\mathfrak{A}_L)\quad\Longleftrightarrow\quad \lambda\in\rho((1+\hat{a}(\lambda))A+Le_\lambda).$$
\end{proof}

\begin{remark}
In the particular case of $a(\cdot)\equiv 0$, we obtain the known result which for $\lambda\in\rho(\mathcal{A}_0)$
$$\lambda\in\rho(\mathfrak{A}_L)\quad\Longleftrightarrow\quad \lambda\in\rho(A+Le_\lambda)$$
which coincides perfectly with the classic problem with delay (see \cite[Lemma 4.1]{BatPia}). On the other hand, for the free-delay integro-differential equation we can identify $\rho(\mathfrak{A}_L)$ with $\rho(\mathcal{A}^{\delta_0})$. Thus for
$\lambda\in {\C}_0\cap\rho(\mathcal{A}_0)$, we retrieve \cite[Proposition 7.25]{EnNag} saying that
$$\lambda\in\rho(\mathcal{A}^{\delta_0})\quad\Longleftrightarrow\quad \lambda\in\rho((1+\hat{a}(\lambda))A).$$
\end{remark}
\section{Intego-differential systems with delays as regular Salamon-Weiss systems}\label{last-sect}
In this section, we keep the same notation as in the previous sections. In addition, set the following Riemann Stieltjes integrals
\begin{align*}
& K\psi=\int^0_{-r}d\mu^K(\theta)\psi(\theta),\quad \mu^{K}:[-r,0]\to\mathcal{L}(U,X)\\ & C\varphi=\int^0_{-r}d\mu^C(\theta)\varphi(\theta),\quad \mu^{C}:[-r,0]\to\mathcal{L}(X,Y),\\ & D\psi=\int^0_{-r}d\mu^D(\theta)\psi(\theta),\quad \mu^{D}:[-r,0]\to\mathcal{L}(U,Y),
\end{align*}
for $\varphi\in W^{1,p}([-r,0],X)$ and $\psi\in W^{1,p}([-r,0],U)$, where $\mu^{K},\,\mu^{C},$ and $\mu^{D}$ are functions of bounded variations assumed continuous on $[-r,0]$ and vanishing at zero.

Now consider the state, input and output delay system
\begin{align}\label{problem-1}
\begin{cases}
\dot{x}(t)=Ax(t)+\displaystyle\int_0^ta(t-s)Ax(s) ds+Lx_t+ Ku_t+f(t),& t\geq 0\\
x(0)=x,\quad x_0=\varphi,\quad u_0=\psi\\
y(t)=Cx_t+Du_t, & t\ge 0,
\end{cases}
\end{align}
for initial data $x\in X,\,\varphi\in L^p([-r,0],X)$ and $\psi\in L^p([-r,0],U),$ where $A,L,a$ and the history function $(t\mapsto x_t)$ are defined in the previous sections, $u\in L^p([-r,+\infty),U$ is the control function, and $t\mapsto u_t$ is the control history function defined by $u_t(\theta)=u(t+\theta)$ for any $t\ge 0$ and $\theta\in [-r,0]$.

In order to state and prove results of this section, let us first recall some notion. On the space $L^p([-r,0],U),$ we define the linear operator
\begin{align*}
Q^U \psi=\psi',\quad D(Q^U)=\{\psi\in W^{1,p}([-r,0],U):\psi(0)=0\}.
\end{align*}
The operator $Q^U$ generates the left shift semigroup $(S^U(t))_{t\ge 0}$ on  $L^p([-r,0],U)$ defined by $(S^U(t)\psi)(\theta)=\psi(t+\theta)$ if $t+\theta\le 0$ and zero if not, for all $t\ge 0$ and $\theta\in [-r,0]$.
In the case of smooth control functions $u,$ we have the following results on the existence of the solution of the problem \eqref{problem-1}.
\begin{theorem}\label{existence-smooth-control}
Let the initial conditions $x\in X,\,\varphi\in L^p([-r,0],U),$ $\psi\in D(Q^U)$, $f\in L^p(\R^+,X)$ and a smooth control function $u\in W^{2,p}([0,+\infty),U)$ with $u(0)=0,$. Then  there exists a unique  solution $x(\cdot)$ of \eqref{problem-1} satisfying $ x_t\in D(L_\Lambda)$ for a.e. $t\ge 0$ and
\begin{align}\label{FVC-input-without}
x(t)=R(t)x+\int_0^tR(t-s)(L_\Lambda x_s+K u_s+f(s)) ds,\quad t\ge 0.
\end{align}
\end{theorem}
\begin{proof}
As in Example \ref{regular-shift-system}, the function $t\mapsto u_t$ is the state of a regular linear system $(Q^U,\beta^U,K)$ with control operator $\beta^U:=(\lambda-Q^U_{-1})e_\lambda$ for $\lambda\in\C$. Thus
\begin{align*}
u_t=S^U(t)\psi+\int^t_0 S^U_{-1}(t-s)(-Q^U_{-1}e_0)u(s)ds,\qquad t\ge 0.
\end{align*}
As $\psi\in D(Q^U)$ and $u\in W^{2,p}([0,+\infty),U)$, an integration by parts shows that $u_t\in W^{1,p}([-r,0],U)$. Thus the function $g(t)=Ku_t$ is well-defined for $t\ge 0$. On the other hand, by \cite[Lemma 3.6]{HMR}, we have $W^{1,p}([-r,0],U)\subset D(K_\Lambda)$ and $g(t)=K_\Lambda u_t$ for a.e. $t\ge 0,$ where $K_\Lambda$ is the Yosida extension of $K$ with respect to $Q^U$. As $t\mapsto K_\Lambda u_t$ is the extended output function of the regular linear system $(Q^U,\beta^U,K)$, we have $g\in L^p_{loc}(\R^+,U)$. If we set
\begin{align*}
\zeta(t)=\left(\begin{smallmatrix} g(t)\\ 0\\0\end{smallmatrix}\right),\quad t\ge 0,
\end{align*}
then the problem \eqref{problem-1} can be reformulated as
\begin{align}\label{N-homo}
\dot{z}(t)=\mathfrak{A}_L z(t)+\zeta(t),\quad z(0)=\left(\begin{smallmatrix} x\\ f\\\varphi\end{smallmatrix}\right),\quad t\ge 0,
\end{align}
where $z(\cdot)$ and $\mathfrak{A}_L$ are given by \eqref{function-z} and \eqref{Generator}, respectively. Now by using the same argument as in the proof of Proposition \ref{sol}, one can see that
the mild solution of the inhomogeneous Cauchy problem \eqref{N-homo} is given by
\begin{align*}
z(t)=\mathbb{T}(t)z(0)+\int^t_0 \mathbb{T}_{-1}(t-s)\mathcal{B} M z(s)ds+\int^t_0 \mathbb{T}(t-s) \left(\mathcal{L}_\Lambda z(s)+\zeta(s)\right)ds,
\end{align*}
for $t\ge 0,$ where $(\mathbb{T}(t))_{t\ge 0}$ is the strongly continuous semigroup given by \eqref{Volterra}. The rest of the proof follows exactly as in the proof of Theorem \ref{new-variation}.
\end{proof}
From our discussion in the proof of Theorem \ref{existence-smooth-control}, we adopt the following definition:
\begin{defi}\label{mild-sol-state-input-delays}
 Let the initial conditions $x\in X,\;\varphi\in L^p([-r,0],X)$ and $\psi\in L^p([-r,0],U)$, and let $f\in L^p(\R^+,X)$.  A mild solution of the integro-differential equation in \eqref{problem-1} is a function $x:[-r,+\infty)\to X$ such that
 \begin{align}\label{vcf-xu}
 x(t)=\begin{cases} R(t)x+\displaystyle \int^t_0 R(t-s)\left(L_\Lambda x_s+K_\Lambda u_s+f(s)\right)ds,& t\ge 0\cr \varphi(t),& a.e.\;t\in [-r,0].
 \end{cases}
 \end{align}
\end{defi}
The following result proves the relationship between smooth input solution and mild solution of the problem \eqref{problem-1}.
\begin{prop}\label{relation-smooth-mild}
  The mild solution of \eqref{problem-1} is limit of a sequence of smooth input solutions of \eqref{problem-1}.
\end{prop}
\begin{proof}
Let $x\in X,\,\varphi\in L^p([-r,0],X),$ $\psi\in L^p([-r,0],U)$ and $u\in L^p_{loc}(\R^+,U)$ such that $u_0=\psi$. Let $x(\cdot)$ the corresponding mild solution of the problem \eqref{problem-1}. We approximate $x,\,\varphi$, $\psi$ and $u$ by sequences $(x^n)_n\subset D(A),\,(\varphi_n)_n\subset D(Q^X),$ $(\psi_n)_n\subset D(Q^U)$ and $u^n\in W^{2,p}_{0,loc}(\R^+,U)$, respectively. For any $t\ge 0,$ define the function
\begin{align*}
u^n_t(\theta)=\begin{cases} u^n(t+\theta),& -t\le \theta\le 0,\cr \psi_n(t+\theta),& -r\le\theta< -t.\end{cases}
\end{align*}
According to Theorem \ref{existence-smooth-control}, the following functions
\begin{align}\label{vcf-xu4}
 x^n(t)=R(t)x+\displaystyle \int^t_0 R(t-s)\left(L_\Lambda x^n_s+K u^n_s+f(s)\right)ds,\quad t\ge 0,\;n\in\N.
 \end{align}
 define smooth input solutions of the problem \eqref{problem-1} with initial conditions $x^n,\varphi_n$ and $\psi_n$.
 Let $\Sigma^U:=(S^U,\Psi^U,\Phi^U,\F^U)$ the regular linear system generated by the triple $(Q^U,\beta^U,K)$. As both functions $(t\mapsto K_\Lambda u_t)$ and $(t\mapsto K u^n_t)$ are output functions of $\Sigma^U$, then
 \begin{equation}\label{Extended-ouput-U}
 \begin{split}
 K_\Lambda u_t&= \left(\Psi^U \psi\right)(t)+\left(\mathbb{F}^U u\right)(t),\cr K u^n_t&= \left(\Psi^U \psi_n\right)(t)+\left(\mathbb{F}^U u^n\right)(t).
 \end{split}
 \end{equation}
Let $\al>0$ be arbitrary and  $t\in [0,\al]$, and set $\gamma:=\max\left(|\mu^K|([-r,0]),|\mu^L|([-r,0])\right)$. By H\"{o}lder inequality,  there exists a constant $c_{\al}>0$ such that
\begin{align*}
\left\|K_\Lambda u_\cdot-Ku^n_\cdot\right\|_{L^p([0,\al],X)}\le \gamma c_{t_0}\left( \|\psi-\psi_n\|_{L^p([-r,0],U)}+\|u-u^n\|_{L^p([0,\al],U)}\right)
\end{align*}
Thus $\|K_\Lambda u_\cdot-Ku^n_\cdot\|_{L^p([0,\al],X)}\to 0$ as $n\to\infty$. Similarly, for any $m,n\in\N,$
\begin{align*}
\left\|L_\Lambda x^n_\cdot-L_\Lambda u^m_\cdot\right\|_{L^p([0,\al],X)}\le \gamma c_{\al}\left( \|\varphi_n-\varphi_m\|_{L^p([-r,0],X)}+\|x^n(\cdot)-x^m(\cdot)\|_{L^p([0,\al],X)}\right).
\end{align*}
On the other hand, by H\"older's inequality, for $m,n\in\N,$
\begin{align*}
&\|R\ast \left((L_\Lambda x^n_\cdot-L_\Lambda x^m_\cdot)+(K_\Lambda u^n_\cdot-K_\Lambda u^m_\cdot)\right)\|_{L^p([0,\al],X)}\cr & \quad \le \al c(\al) \left( \|\varphi_n-\varphi_m\|_p+ \| x^n(\cdot)-x^m(\cdot)\|_{L^p([-r,0],X)}+ \|K_\Lambda u^n_\cdot-K_\Lambda u^m_\cdot\|_{L^p([0,\al],X)} \right)
\end{align*}
for a constant $c(\al)>0$. We choose $\al>0$ such that $\al c(\al)<\frac{1}{2}$. For a such $\al,$ we have
\begin{align*}
&\| x^n(\cdot)-x^m(\cdot)\|_{L^p([-r,0],X)}\cr & \qquad\le  2\kappa(\al) \left( \|x^n-x^m\|+ \|\varphi_n-\varphi_m\|_p+ \|K_\Lambda u^n_\cdot-K_\Lambda u^m_\cdot\|_{L^p([0,\al],X)}\right).
\end{align*}
We deduce that $(x^n(\cdot))_n$ is a Cauchy sequence in $L^p([0,\al],X)$, so it converges to function $x(\cdot)\in L^p([0,\al],X)$. We select
\begin{align*}
x_t(\theta)=\begin{cases} x(t+\theta),& -t\le \theta\le 0,\cr \varphi(t+\theta),& -r\le\theta< -t.\end{cases}
\end{align*}
As above,  $L_\Lambda x^n_\cdot$ converges to $L_\Lambda x_\cdot$ in $L^p([0,\al],X)$. By passing to limit in \eqref{vcf-xu4}, we obtain
\begin{align*}
 x(t)=R(t)x+\displaystyle \int^t_0 R(t-s)\left(L_\Lambda x_s+K_\Lambda u_s+f(s)\right)ds,\quad t\ge 0.
 \end{align*}
 this ends the proof.
\end{proof}
In the rest of this section we will show that the delay system \eqref{problem-1} is equivalent to a regular distributed linear system in the Salamon-Weiss sense (see Section \ref{sec:2} for the definition of such systems). In fact, as shown in the proof of Theorem \ref{existence-smooth-control}, the problem \eqref{problem-1} is reformulated as the following system
\begin{equation}\label{JJJ}
\dot{z}(t)=\mathfrak{A}_L z(t)+\left(\begin{smallmatrix} K_\Lambda u_t\\0\\0\end{smallmatrix}\right),\quad z(0)=\left(\begin{smallmatrix} x\\ f\\\varphi\end{smallmatrix}\right),\quad t\ge 0.
\end{equation}
In particular, the mild solution of the equation \eqref{JJJ} is given by
\begin{equation}\label{z-solution}
z(t)=\mathfrak{T}_L(t)\left(\begin{smallmatrix} x\\ f\\\varphi\end{smallmatrix}\right)+\int^t_0\mathfrak{T}_L(t-s)\left(\begin{smallmatrix} K_\Lambda u_s\\0\\0\end{smallmatrix}\right)ds.
\end{equation}
Let us introduce the a new product state space
\begin{equation*}
\tilde{\mathcal{X}}=\mathcal{X}\times L^p([-r,0],U),
\end{equation*}
where $\mathcal{X}$ is the product space defined in \eqref{space-calX}. Moreover, we select a new state function
\begin{equation*}
w:t\in [0,+\infty)\mapsto w(t)=\begin{pmatrix} z(t)\\u_t\end{pmatrix}\in \tilde{\mathcal{X}}.
\end{equation*}
The main result of this section is the following:
\begin{theorem}
  There exist a generator $\tilde{\mathfrak{A}}:D(\tilde{\mathfrak{A}})\subset \tilde{\mathcal{X}}\to \tilde{\mathcal{X}} $ of a strongly continuous semigroup on $\tilde{\mathcal{X}},$ an admissible control operator $\tilde{\mathcal{B}}\in\mathcal{L}(U,\tilde{\mathcal{X}}_{-1})$ for  $\tilde{\mathfrak{A}}$ and an admissible observation operator $\tilde{\mathcal{C}}\in\mathcal{L}(D(\tilde{\mathfrak{A}}),Y)$ for $\tilde{\mathfrak{A}}$ such that the triple $(\tilde{\mathfrak{A}},\tilde{\mathcal{B}},\tilde{\mathcal{C}})$ is regular on $\tilde{\mathcal{X}},U,Y$ and the integro-differential equation with state and input delays \eqref{problem-1} is reformulated as the following input-output distributed linear system
  \begin{equation*}
  \begin{cases}
  \dot{w}(t)=\tilde{\mathfrak{A}}w(t)+\tilde{\mathcal{B}} u(t),\quad w(0)=\left(\begin{smallmatrix} x\\ f\\\varphi\\\psi\end{smallmatrix}\right),& t\ge 0,\\
  y(t)= \tilde{\mathcal{C}} w(t),& t\ge 0.
  \end{cases}
  \end{equation*}
\end{theorem}
\begin{proof}
By using \eqref{Extended-ouput-U}, \eqref{z-solution} and the notation
\begin{align*}
& \mathscr{M}(t)\psi:=\int^t_0\mathfrak{T}_L(t-s)\left(\begin{smallmatrix} (\Psi^U \psi)(s)\\0\\0\end{smallmatrix}\right)ds,\cr & \mathscr{N}_t u:=\int^t_0\mathfrak{T}_L(t-s)\left(\begin{smallmatrix} (\F^U u)(s)\\0\\0\end{smallmatrix}\right)ds
\end{align*}
for $\psi\in L^p([-r,0],U)$ and $u\in L^p(\R^+,U)$ with $u_0=\psi,$ we can also write
\begin{equation*}
z(t)=\mathfrak{T}_L(t)\left(\begin{smallmatrix} x\\ f\\\varphi\end{smallmatrix}\right)+\mathscr{M}(t)\psi+\mathscr{N}_t u.
\end{equation*}
Now, the fact that $u_t=S^U(t)\psi+\Phi^U_t u,$ where $(\Phi^U_t)_{t\ge 0}$ is the family of control maps associated with the regular triple $(Q^U,\beta^U,K)$, implies that
\begin{equation*}
w(t)=\tilde{\mathfrak{T}}(t)\left(\begin{smallmatrix} x\\ f\\\varphi\\\psi\end{smallmatrix}\right)+\tilde{\Phi}_t u,
\end{equation*}
where
\begin{equation*}
\tilde{\mathfrak{T}}(t):=\left(\begin{array}{c|c}
                                   \mathfrak{T}_{L}(t) & \mathscr{M}(t) \\
                                   \hline
                                   \begin{matrix}
                                     0 & 0 & 0
                                   \end{matrix} & S^U(t)
                                 \end{array}\right),\quad\text{and}\quad \tilde{\Phi}_t u:=\begin{pmatrix}\mathscr{N}_t u\\ \Phi^U_t u\end{pmatrix}.
\end{equation*}
A similar argument as in the proof of \cite[Theorem 3.1]{Ha-Id-IMA} shows that the operators family $(\tilde{\mathfrak{T}}(t))_{t\ge 0}$ define a strongly continuous semigroup on $\tilde{\mathcal{X}}$ with generator
\begin{equation*}
\tilde{\mathfrak{A}} =\left(\begin{array}{c|c}
                          \mathfrak{A}_{L} & \begin{matrix}
                                               K \\
                                               0 \\
                                               0
                                             \end{matrix} \\
                                             \hline
                          \begin{matrix}
                            0 & 0 & 0
                          \end{matrix} & Q^U
                        \end{array}\right),\quad
  D( \tilde{\mathfrak{A}}) =D( \mathfrak{A}_{L})\times D(Q^U).
\end{equation*}
We know that (see Example \ref{regular-shift-system})
\begin{align*}
\|\Phi^U u\|_{L^p([-r,0],U)}\le \|u\|_{L^p([0,t],U)}\quad\text{and}\quad \|\mathbb{F}^Uu\|_{L^p([0,t],U)}\le |\mu^K|([-r,0])\|u\|_{L^p([0,t],U)}
\end{align*}
for any $t\ge 0$ and $u\in L^p(\R^+,U)$. Then by H\"older inequality, we obtain
\begin{align*}
\|\tilde{\Phi}_t u\|_{\tilde{\mathcal{X}}}&= \|\mathscr{N}_tu\|_{\mathcal{X}}+\|\Phi^U u\|_{L^p([-r,0],U)}\cr & \le (c  |\mu^K|([-r,0])+1 )\|u\|_{L^p([0,t],U)}
\end{align*}
for some constant $c:=c(t)>0$. In addition, the fact that
\begin{align*}
\Phi^U_{t+s}u=S^U(t)\Phi^U_s (u_{|[0,s]})+\Phi^U_t u(\cdot+s),\\
(\mathbb{F}^U u)(t+s)=(\Psi^U\Phi_s u)(t)+(\mathbb{F}^U u(\cdot+s))(t),
\end{align*}
implies that
\begin{align*}
\tilde{\Phi}_{t+s}u=\tilde{\mathfrak{T}}(t)\tilde{\Phi}_{s}u+\tilde{\Phi}_{t}u(\cdot+s).
\end{align*}
According to \cite{Weiss-control-89}, this functional equation shows that there exists an admissible control operator $\tilde{\mathcal{B}}\in\mathcal{L}(U,\tilde{\mathcal{X}})$ for $\tilde{\mathfrak{A}}$ such that
\begin{equation*}
\tilde{\Phi}_t u=\int^t_0 \tilde{\mathfrak{T}}_{-1}(t-s)\tilde{\mathcal{B}} u(s)ds
\end{equation*}
for any $t\ge 0$ and $u\in L^p(\mathbb{R}^+,U)$. Thus the function $w(\cdot)$ is the solution of the control problem
\begin{align*}
 \dot{w}(t)=\tilde{\mathfrak{A}}w(t)+\tilde{\mathcal{B}} u(t),\quad w(0)=\left(\begin{smallmatrix} x\\ f\\\varphi\\\psi\end{smallmatrix}\right),\quad t\ge 0
\end{align*}
We select
\begin{align*}
\tilde{\mathcal{C}}:=\begin{pmatrix} 0&0&C&D\end{pmatrix}:D(\tilde{\mathfrak{A}})\to Y.
\end{align*}
Then the out-put of the system \eqref{problem-1} satisfies
\begin{align*}
y(t)=\tilde{\mathcal{C}} w(t),\qquad t\ge 0.
\end{align*}
First, we prove that $\tilde{\mathcal{C}}$ is an admissible observation operator for $\tilde{\mathfrak{A}}$. According to the expression the semigroup $(\tilde{\mathfrak{T}}(t))_{t\ge 0},$ the fact that the operator $D\in\mathcal{L}(D(Q^U),Y)$ is admissible observation for $Q^U,$ and \cite[Proposition 3.3]{Hadd-SF} it suffices to show that the operator $\mathcal{C}=(0~0~C)\in\mathcal{L}(D(\mathfrak{A}),Y)$ is an admissible observation operator for $\mathfrak{A}$. To this end, as $\mathfrak{A}_L=\mathfrak{A}+\mathcal{L}$ and $\mathcal{L}$ is an admissible observation operator for $\mathfrak{A},$ (see the proof of Theorem \ref{generation-big-A}), by \cite{Ha-Id-SCL} it suffices to prove that $\mathcal{C}$ is an admissible observation operator for $\mathfrak{A}$. This is obvious due to the expression of the semigroup $(\mathbb{T}(t))_{t\ge 0}$ given in \eqref{Volterra} and the fact that the operator $C\in\mathcal{L}(D(Q),X)$ is an admissible observation operator for the left shift semigroup $(S(t))_{t\ge 0}$ on $L^p([-r,0],X)$.\\ Second, we prove that the triple $(\tilde{\mathfrak{A}},\tilde{\mathcal{B}},\tilde{\mathcal{C}})$ is regular. To that purpose, we need  first to compute the Yosida extension of $\tilde{\mathcal{C}}$ with respect to $\tilde{\mathfrak{A}}$. In fact, by taking the Laplace transform of $\tilde{\mathfrak{T}}(t)$ we obtain
\begin{align*}
R(\lambda,\tilde{\mathfrak{A}})=\left(\begin{array}{c|c}
                                   R(\lambda,\mathfrak{A}_{L}) & R(\lambda,\mathfrak{A}_{L})\left(\begin{smallmatrix}KR(\lambda,Q^U)\\0\\0\end{smallmatrix}\right) \\
                                   \hline
                                   \begin{matrix}
                                     0 & 0 & 0
                                   \end{matrix} & R(\lambda,Q^U)
                                 \end{array}\right),\qquad \lambda\in \rho(\mathfrak{A}_{L}).
\end{align*}
Let $(x~f~\varphi~\psi)^\top\in\tilde{\mathcal{X}}$. For $\lambda>0$ sufficiently large we have
\begin{align*}
\tilde{\mathcal{C}}\lambda R(\lambda,\tilde{\mathfrak{A}})\left(\begin{smallmatrix}x\\f\\\varphi\\\psi\end{smallmatrix}\right)= \mathcal{C}\lambda R(\lambda,\mathfrak{A}_{L})\left(\begin{smallmatrix}x\\f\\\varphi\end{smallmatrix}\right)+\mathcal{C} R(\lambda,\mathfrak{A}_{L})\left(\begin{smallmatrix}K\lambda R(\lambda,Q^U)\psi\\0\\0\end{smallmatrix}\right)+D\lambda R(\lambda,Q^U)\psi.
\end{align*}
Observe that if $\psi\in D(K_\Lambda)$, then
\begin{align*}
\left\|\mathcal{C} R(\lambda,\mathfrak{A}_{L})\left(\begin{smallmatrix}K\lambda R(\lambda,Q^U)\psi\\0\\0\end{smallmatrix}\right)\right\|\le \left\|\mathcal{C} R(\lambda,\mathfrak{A}_{L})\right\|\left(\|K\lambda R(\lambda,Q^U)\psi -K_\Lambda \psi \|+\|K_\Lambda \psi\|\right).
\end{align*}
As $\mathcal{C}$ is an admissible observation operator for $\mathfrak{A}_{L},$ it follows that $\|\mathcal{C} R(\lambda,\mathfrak{A}_{L})\|$ goes to $0$ when $\lambda\to+\infty$. Then
\begin{align*}
\lim_{\lambda\to+\infty} \left\|\mathcal{C} R(\lambda,\mathfrak{A}_{L})\left(\begin{smallmatrix}K\lambda R(\lambda,Q^U)\psi\\0\\0\end{smallmatrix}\right)\right\|=0,\qquad \forall \psi\in D(K_\Lambda).
\end{align*}
Then we have
\begin{align}\label{inclusion}
\Omega:=D(\mathcal{C}_\Lambda)\times (D(K_\Lambda)\cap D(D_\Lambda))\subset D(\tilde{\mathcal{C}}_\Lambda)\quad\text{and}\quad (\tilde{\mathcal{C}}_\Lambda)_{|\Omega}=\begin{pmatrix} \mathcal{C}_\Lambda & D_\Lambda\end{pmatrix}.
\end{align}
By \cite[Proposition 3.3]{Hadd-SF}, for $u\in L^p_{loc}(\mathbb{R}^+,U)$ we have
\begin{align*}
\mathscr{N}_t u\in D(\mathcal{C}_\Lambda)\quad\text{and}\quad \|\mathcal{C}_\Lambda\mathscr{N}_\cdot u\|_{L^p([0,\tau],Y)}&\le c \tau^{\frac{1}{q}} \|\F^Uu\|_{L^p([0,\tau],X)}\cr &\le c \tau^{\frac{1}{q}} |\mu^K|([-\tau,0]) \|u\|_{L^p([0,\tau],U)},
\end{align*}
for a.e. $t\ge 0$ where $c>0$ is a constant and $\frac{1}{p}+\frac{1}{q}=1$. On the other hand, the triple $(Q^U,\beta^U,K)$ and $(Q^U,\beta^U,D)$ are regular and having the same control maps. Then $\Phi^U_t u\in  D(K_\Lambda)\cap D(D_\Lambda)$ for a.e. $t\ge 0$. Moreover, the extended input-output operator $\F^{D}$ of $(Q^U,\beta^U,D)$ is given by $(\F^{D}u)(t)=D_\Lambda \Phi^U_t u$ for almost every $t\ge 0,$ and
\begin{align*}
\|\F^{D}u\|_{L^p([0,\tau],Y)}\le |\mu^D|([-\tau,0]) \|u\|_{L^p([0,\tau],U)}.
\end{align*}
Now according to \eqref{inclusion}, we have
\begin{align*}
\tilde{\Phi}_t u\in \Omega \subset D(\tilde{\mathcal{C}}_\Lambda)\quad\text{and}\quad (\tilde{\F}u):=\tilde{\mathcal{C}}_\Lambda \tilde{\Phi}_t u= \mathcal{C}_\Lambda \mathscr{N}_tu+(\F^{D}u)(t)
\end{align*}
for a.e. $t\ge 0$ and all $u\in L^p_{loc}(\mathbb{R}^+,U)$. Moreover
\begin{align}\label{last-estimate}
\|\tilde{\F}u\|_{L^p([0,\tau],Y)}\le c(\tau) \|u\|_{L^p([0,\tau],U)},\qquad \forall u\in L^p([0,\tau],U),
\end{align}
where
\begin{align*}
c(\tau):=2^{p-1}\left(c \tau^{\frac{1}{p}} |\mu^K|([-\tau,0])+|\mu^D|([-\tau,0])\right)\underset{\tau\to 0}{\longrightarrow} 0.
\end{align*}
Hence the triple $(\tilde{\mathfrak{A}},\tilde{\mathcal{B}},\tilde{\mathcal{C}})$ is well-posed in $\tilde{\mathcal{X}},U,Y$. Let $b\in U$ be a fixed control and set $u_b(t)=b$ for any $t\ge 0$. By using H\"older inequality and the estimate \eqref{last-estimate}, we obtain
\begin{align*}
\left\|\frac{1}{\tau}\int^{\tau}_0 (\tilde{\F}u_b)(t)dt\right\|\le c(\tau) \|b\|_U.
\end{align*}
This implies that
\begin{align*}
\lim_{\tau\to 0}\frac{1}{\tau}\int^{\tau}_0 (\tilde{\F}u_b)(t)dt=0.
\end{align*}
Thus the triple $(\tilde{\mathfrak{A}},\tilde{\mathcal{B}},\tilde{\mathcal{C}})$ is regular in $\tilde{\mathcal{X}},U,Y$.

\end{proof}

\end{document}